\documentclass[preprint,12pt]{article}
\usepackage{graphicx}

\usepackage{amssymb}
\usepackage{amsmath}

\usepackage{xcolor}

\begin{document}




\title{Model discovery on the fly using continuous data assimilation}
\author{J. Newey, J. P. Whitehead, and E. Carlson}




\maketitle

\begin{abstract}
We review an algorithm developed for parameter estimation within the Continuous Data Assimilation (CDA) approach.
We present an alternative derivation for the algorithm presented in a paper by Carlson, Hudson, and Larios (CHL) \cite{CHL_2020}.
This derivation relies on the same assumptions as the previous derivation but frames the problem as a finite dimensional root-finding problem.
Within the approach we develop, the algorithm developed in \cite{CHL_2020} is simply a realization of Newton's method. 
We then consider implementing other derivative based optimization algorithms; we show that the Levenberg Maqrquardt algorithm has similar performance to the CHL algorithm in the single parameter estimation case and generalizes much better to fitting multiple parameters.
We then implement these methods in three example systems: the Lorenz '63 model, the two-layer Lorenz '96 model, and the Kuramoto-Sivashinsky equation.
\end{abstract}









\section{Introduction}

Fully data driven methods such as deep neural networks have recently advanced simulation methods for physical models (see \cite{bi2023accurate} for example).
However even modified versions of these methods that incorporate the physical equations of motion into the requisite loss function \cite{raissi2019physics} can be difficult to interpret.  At the same time, traditional physical models that are derived from first principles frequently do not match the observable data well, due to a mismatch between the parameters of the model and their physically realizable values. Alternative hybrid approaches have been developed that use existing data to identify the form of the governing equations for a given system (see \cite{brunton2016discovering} for example), or use incoming data to augment the form of an existing model (see \cite{chen2023data,mojgani2024interpretable,guan2024online,jakhar2024learning}), and other studies have investigated the sensitivity and utility of parameters in existing models using observable data to identify the `true' or at least optimal parameter that best fits the data and underlying physics simultaneously (see \cite{burger2014maximum,jiang2022embed,rumbell2023novel,pilosov2023parameter,miguez2023sequential}).  Most of these approaches rely on statistical arguments that require ensembles of simulations which quickly become computationally prohibitive.  In addition the equation discovery and/or parameter estimation/update is performed posteriori, i.e. after the data has been collected and processed.

In this article we present a new framework for algorithms that estimate unknown parameters in a dynamical system using partial observations of the state, which circumvent both of these issues: 1) only a single real-time simulation is required, and 2) the parameters of the model are updated as the data is recorded/observed. The framework developed here is built on the continuous data assimilation setting first introduced in \cite{azouani2014continuous} and further explored in  several resources including \cite{gesho2016computational,bessaih2015continuous,farhat2015continuous,albanez2016continuous,farhat2016abridged,farhat2016data,farhat2016charney,foias2016discrete,biswas2017higher,jolly2017data,jolly2017determining,albanez2018continuous,blocher2018data,farhat2018assimilation,mondaini2018uniform,celik2019spectral,jolly2019continuous,clark2020synchronization,farhat2020data,ibdah2020fully,biswas2021data,cao2021algebraic,cao2022continuous,franz2022bleeps,celik2023data,jolly2023data,larios2023second,larios2023application,biswas2024unified,diegel2024analysis,larios2024continuous,you2024continuous,carlson2024determining}.  Recently this setting has been modified to develop algorithms that identify unknown parameters in the dynamical system \cite{clark2018inferring,CHL_2020,carlson2021dynamically,martinez2022convergence,martinez2022reconstruction,biswas2023determining,albanez2024parameter,ccibik2024adaptive,farhat2024identifying,lu2024continuous,martinez2024relaxation}.  We demonstrate that under some reasonable assumptions, the algorithm developed in \cite{CHL_2020} reduces to Newton's method for rootfinding acting on a loss function defined as the $L^2$ norm of the observable error in the state.
Using the general optimization framework we develop here, we additionally show that by applying the Gauss Newton method for optimization we can identify a parameter estimation scheme which better generalizes to a multi-dimensional parameter space and un-parameterized model error.
The key result is that because the optimal relaxation parameter $\mu$ specified first in \cite{azouani2014continuous} is very large, it allows us to do perturbation theory on the sensitivity equations giving a simple, on-line approximation of the gradient of the underlying loss function.


\section{Error Function and Sensitivity Equations}

\subsection{Problem Statement}

We begin by considering a modification of the continuous data assimilation (CDA) approach where the model is no longer exactly known. 
We make the modification that the simulated system has some model error which takes the form of a finite number of unknown parameters.
To denote this,  write $F(u;\gamma): H \times \mathbb{R}^n \rightarrow H$ where in this article $H=\mathbb{R}^N$ for some $N\neq n$, but in general $H$ may be an appropriate functional space.
The CDA motivated system we consider is given by
\begin{align}
    &\Dot{u} + F(u; \mathbf{\gamma}) = 0 \label{eq:original_system}\\
    &\Dot{v} + F(v; c)  + \mu I_h(v-u) = 0,\label{eq:nudge_system}
\end{align}
where $\gamma = (\gamma_i) \in \mathbb{R}^n$ is some vector of ``true'' parameters which are assumed unknown, and $c  = (c_i)\in \mathbb{R}^n$ is a vector of approximate parameters.
We use the notation $\Dot{u} = \frac{d u }{d t}$.
The ``true'' state of the system is $u(t)$ and the data assimilated state is $v(t)$.  When $c=\gamma$, i.e. the model is fully known to the user, \cite{azouani2014continuous} demonstrate that under suitable restrictions on the system and the observation operator $I_h$ (denoting the portion of the state that is available to the user at any given time $t$), $v(t)\rightarrow u(t)$ as $t\rightarrow \infty$, in appropriate norms.  The problem when $c\neq \gamma$ was originally considered in \cite{farhat2020data} and more particularly in \cite{CHL_2020} with further extensions in \cite{pachev2022concurrent} and \cite{martinez2024relaxation}.
In the following we will develop a general parameter identification/estimation approach for which \cite{CHL_2020} is a special case.

Our goal is to find the ideal parameter values that minimize the error functional
\begin{align}
    E(c) = \frac{1}{2}\|I_h(v(c)-u)\|^2_H \bigg|_{t=T>>1}.
\end{align} 
The derivation below is not particular to the choice of norm in defining the error functional $E(c)$; in fact an almost identical derivation will apply to other cost functionals so long as they satisfy some basic assumptions such as convexity and Fr\'{e}chet differentiability.  The approach outlined below is also applicable when the actual form of the ``true'' system is unknown that is when not only are the parameters $\gamma$ unknown, but the functional form of $F(u;\gamma)$ is unknown as well.  This is demonstrated in Section \ref{sec:optvsroot}. 

\subsection{Assumptions}
We make several key assumptions in the following derivation which are outlined here. 
These same assumptions generally underlie many similar parameter estimation approaches such as \cite{carlson2021dynamically}.
These assumptions were originally motivated by numerical experimentation, but have been justified more rigorously in several cases \cite{martinez2022convergence,CHL_2020}.
The key idea is that given these assumptions, the problem reduces to a finite dimensional nonlinear optimization problem.

\begin{figure}
    \centering
    \includegraphics[width=0.9\textwidth]
    {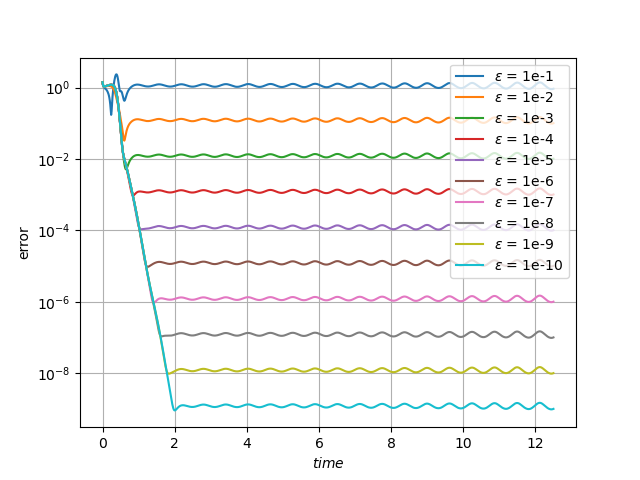}
    \caption{Time independence of long time error for the Lorenz '63 system while varying across an approximate parameter $c_1$,
    where $\epsilon = \frac{\Delta c_1}{c_1}$.
    We see that the nudging converges up to some long term error proportional to the parameter error.
    }
    \label{fig:TimeInd}
\end{figure}

\subsubsection{Assumption 1: Time independence of long time error}
Consistent with observations first noted in \cite{farhat2020data} and \cite{CHL_2020} we anticipate that when $c\neq \gamma$ then $v(t)$ will converge close to $u(t)$ up to an error term that is proportional to the parameter error.  This is illustrated for the Lorenz '63 system as shown in Figure \ref{fig:TimeInd}, and motivates 
our assumption that for some $T \gg 1$ we have
\begin{align}
    \frac{\partial}{\partial t} E(c)
    \bigg|_{t=T} 
    \approx 0.
\end{align}
This assumes that the system has been nudged for sufficient time such that the error has settled to a steady state value.  This motivates the choice of notation $E(c)$, which we assume is time-independent for sufficiently large times.
In practice the system will tend to relax and then fluctuate around a set value, i.e. the time-independence assumption may not be always realized in practice (see Figure \ref{fig:TimeInd} for example). 
A more robust modification of the following that accounts for these fluctuations would include time integrals that average away such fluctuations.
Such a modification of parameter estimation algorithms is mentioned in \cite{CHL_2020, martinez2022convergence}.

\subsubsection{Assumption 2: Independence from initial conditions}
\begin{figure}
    \centering
    \includegraphics[width=0.9\textwidth]
    {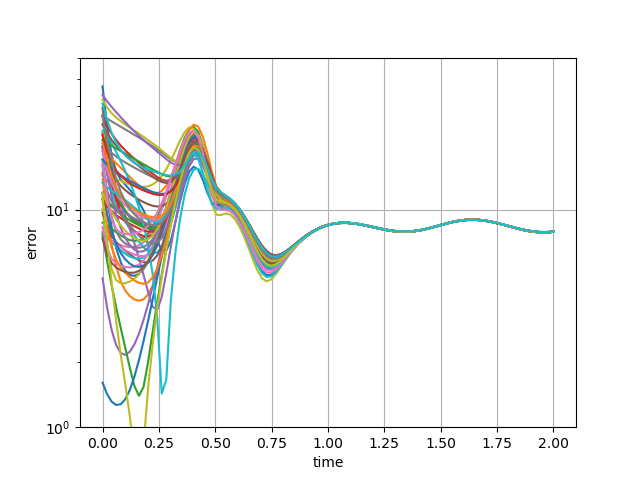}
    \caption{Error between the true state and the assimilated state for 50 randomly chosen initial conditions for the assimilated state when the parameter error is $\epsilon = 0.5$ for the Lorenz '63 model.  Note that despite the differences in the initial conditions, after a sufficiently long time (approximately $1.0$ here) the error collapses to a value dictated by the parameter/model error.
    }
    \label{fig:InitCond}
\end{figure}

We make the key assumption that E(c) is independent of the initial condition.  While not entirely true for any given $t =T$, the large time behavior ($t \to \infty$) of the error is independent of the initial conditions taken within the absorbing ball, whether or not the parameters are correct (see, e.g., \cite{farhat2020data, CHL_2020, biswas2023determining, albanez2024parameter, martinez2024relaxation}.  Thus, for any initial condition within the absorbing ball, the methods employed here will converge, as we demonstrate computationally for the Lorenz '63 system in Figure \ref{fig:InitCond}

\subsubsection{Assumption 3: Derivatives of  E(c) are continuous}
Our goal is to implement gradient based optimization and root finding methods on $E(c)$.
This requires estimating the derivatives of $E(c)$.  
Hence we will assume that $E(c)$ is continuously differentiable in the parameters $c$.  In addition, to make the differentiation more practical we assume that the observation operator $I_h$ is a linear projection operator and independent of the unknown parameters $c$.
We make this assumption mainly for convenience, though this approach could be generalized to the case where observation function is nonlinear.

\subsection{Sensitivity equations}
For a specific parameter $c_i$ we write out the derivative of $E$ with respect to $c_i$ as function of $w_i := \frac{\partial v}{\partial c_i}$,
\begin{align}
    \frac{\partial E}{\partial c_i}
     = 
     \frac{1}{2}
     \frac{\partial}{\partial c_i}
     \| I_h(v-u)\|_H^2
     \bigg|_{t=T}
     = \langle v-u, I_h w_i \rangle_H \bigg|_{t=T}. 
\end{align}
We find the time evolution of $w_i$ by implicitly taking the derivative of \eqref{eq:nudge_system} with respect to $c_i$.
From assumption 3 above, we have implicitly assumed differentiability of $F$ in both $v$ and $c$.
If our space of solutions $H$ is infinite dimensional then a rigorous analysis would be necessary to justify Fr\'{e}chet differentiability of $F$ with respect to the solution $v$.
We move forward assuming that the formal manipulations performed here are justified.
This leads us to the sensitivity equations (see \cite{dickinson1976sensitivity}) for the nudged system \eqref{eq:nudge_system},
\begin{align}\label{eq:sensitivity1}
    \Dot{w}_i + D F(v;c)
     w_i
    + 
    F_{c_i}(v; c)  + \mu I_h w_i
     = 0.
\end{align} 
Here, $DF$ is the Jacobian of $F$  in $v$, and $F_{c_i}$ is the partial derivative of $F$ with respect to $c_i$ with $v$ held constant. 
We will assume by default that the initial conditions do not depend on the parameters $c$ so that the initial condition is given by:
\begin{align}
    w_i |_{t= 0}
    = 0.
\end{align}

Note that \eqref{eq:sensitivity1} is linear in the $w_i$, although coupled to the (likely nonlinear) \eqref{eq:nudge_system} through the nudged state $v(t)$. 
The sensitivity equations are also greatly simplified by having a zero initial condition.
In addition, note that $\mu$ incorporates a decay term into the sensitivity equations, which doesn't allow the solution of \eqref{eq:sensitivity1} to grow in time, i.e. the sensitivities of the solution $v(t)$ to the parameters $c$ can be controlled by the size of $\mu$, at least for large times.

\subsection{Asymptotics\label{sec:asym}}
Making use of the sensitivity equations in a gradient
based optimization routine to identify a minimum of $E(c)$ is practical for systems of ordinary differential equations (ODE) such as the Lorenz '63 system, but for partial differential equations (PDE) where a single simulation of \eqref{eq:nudge_system} is computationally prohibitive, the added cost of the sensitivity equations, makes such a parameter estimation algorithm implausible. 
The number of sensitivity equations grows with the number of parameters, and so this makes direct simulation computationally prohibitive in systems with a large number of parameters. 
This problem could be reduced by using an adjoint method, but we will not explore this here. 
Directly simulating the sensitivity equation for parameter learning has been done in systems biology \cite{raue2013lessons,frohlich2017parameter}, however not concurrently with the CDA approach as done here.

Remarkably, relatively simple asymptotic approximations of the sensitivity equations yield approximate gradients that imitate the true gradient extremely well.
In order to keep the derivation general, the calculation in this section is largely formal.
For finite dimensional systems, making these asymptotic approximations is a straightforward application of some standard techniques, whereas in the infinite dimensional, i.e. partial differential equation (PDE), case the analysis will be significantly more involved.

Since the convergence theory requires that $\mu$ is sufficiently large to theoretically guarantee convergence, we consider the asymptotics for the limit as $\mu \to \infty$. 
We will approach this by first rescaling the time variable as $\tau = \mu t$. 
Thus we write
\begin{align}
    w_i(t) = W_i(\tau).
\end{align}
By the chain rule we see that we have
\begin{align}
    \Dot{w}_i
    = \mu \Dot{W} .
\end{align}
Plugging this into the sensitivity equations, we get
\begin{align}
    \mu \Dot{W}_i
    + D F W_i
    + \mu I_h W_i
    = -  F_{c_i}.
\end{align}
If we let $W = W^0 + \mu^{-1} W^1 + \mu^{-2} W^2+...$  and follow standard perturbation theory, we find that $I_h W_i^0 = 0$, and the next term $W_i^1$ satisfies:
\begin{align}
    &\Dot{W}^1_i
    + I_h (W^1_i)
    = -  F_{c_i}(v), \notag\\
    &W^1_i(0) = 0.
\end{align}
The solution, in the $I_h$ (observable) subspace, is given by,
\begin{align}
    I_h W^1_i(\tau)
    = - \int_0^\tau I_h F_{c_i}
    \left(v\left(\tau'/\mu\right)\right) e^{\tau'-\tau} 
    d\tau'.
\end{align}

This gives us the following first order asymptotic approximation for $w_i$:
\begin{align}
    I_hw_i(t) \sim -\frac{1}{\mu} \int_0^{\mu t} I_h F_{c_i}(v(\tau'/\mu)) 
    e^{\tau' - \mu t} d\tau'.
\end{align}
If we make the change of variables $s =  - \tau' /\mu + t $, and directly apply Watson's Lemma \cite{miller2006applied} then we can identify the leading order approximation: 
\begin{align}\label{eq:asymptotic1}
    I_hw_i(t)
    \sim -\frac{1}{\mu} I_hF_{c_i}(v(t)).
\end{align}
We have also considered the second derivative sensitivity equations which yield higher order optimization methods, however we do not include these results here as such methods do not seem to lead to noticeable improvements upon the leading order approximation.

\section{Parameter estimation}
When there is a unique exact parameter $\gamma$ for the system \eqref{eq:original_system} then $E(c) = \frac{1}{2}\|I_h(v(c)-u)\|^2 \big|_{t = T \gg 1}$ as a function of $c$ has a unique root and global minima at $c = \mathbf{\gamma}$.
Hence, taking advantage of the asymptotics from the previous section, we can use any gradient based root-finding or optimization routine to develop a parameter update algorithm seeking to find the root and/or minimize $E(c)$.
We follow the ``relax then punch'' approach where we allow the error to relax between parameter update steps \cite{CHL_2020,carlson2021dynamically,martinez2022convergence,farhat2024identifying,martinez2024relaxation}.
We require a gap between parameter updates so that Assumptions 1 and 2 are adequately satisfied.
This amounts to the assumption that $\Delta t = t_{k+1} - t_k$ is ``large'' (where $t_k$ is our parameter update time). 
In practice, the system usually relaxes very quickly, so the interval between parameter updates does not need to be especially large.
In fact, we note from prior work on CDA that $\| v - u\| \sim e^{-\mu t}$ so that a  characteristic relaxation time scale is $\frac{1}{\mu}$. Thus we choose $\Delta t$ much larger than this characteristic time scale.
In practice because systems tend to fluctuate around a minimum rather than completely relax, there seems to be a sweet spot for the value of $\Delta t$ where the error between the systems has just relaxed but not started to fluctuate.
A good rule of thumb seems to be choosing $\Delta t  >  \frac{10}{\mu}$.


As a technical note, because the initial conditions for each simulation depend on the final conditions of the previous run, the initial conditions for the sensitivities will not in general be zero and will be hard to estimate a priori (They can be computed numerically of course). 
This isn't as problematic as it first appears, as we could restart after each parameter update with random initial conditions and (assuming the validity of assumption 2) the parameter update algorithm will still work.
This is observed numerically in all three of the systems considered in this paper (not shown).
Moreover because the sensitivity equations are linear (taking $v$ as given) if we have initial conditions $w_i(0)$ for the sensitivities we write, 
\begin{align}
    w_i = w_{i,1} + w_{i,2}
\end{align}
where, 
\begin{align}
    \dot{w}_{i,1}  + DF(v)w_{i,1} +\mu I_h w_{i,1} &= 0 \text{ and } w_{i,1}(0) = w_i(0)\\
    \dot{w}_{i,2}  + DF(v) w_{i,2}+\mu I_h w_{i,2} &= -F_{c_i}(v) \text{ and } w_{i,2}(0) = 0
\end{align}
So long as the Jacobian $DF$ is bounded in time and $\mu$ is sufficiently large, $I_h w_{i,1}$ exponentially decays in time, and hence because of the large time assumption the initial sensitivity dependent terms can be ignored.

In all of the methods demonstrated below, we require the computation or at least approximation of the sensitivities $w_i$.  To make a robust comparison, we will include some calculations where we have directly simulated (DS) the sensitivity equations \eqref{eq:sensitivity1}, and other cases where we have approximated the sensitivities via the first order asymptotic representation in \eqref{eq:asymptotic1} which we refer to as ``on the fly'' (OTF) methods. As we demonstrate below, the OTF method has similar performance to DS method at a much lower computational cost.

\subsection{Root finding methods, and the re-derivation of the algorithm proposed in \cite{CHL_2020}}
The accelerated Newton's method for finding a root of order $m$ of the scalar valued function $E(c)$ is given in \cite{frontini2003modified}:
\begin{align}
    c^{(k+1)} = c^{(k)} - m\frac{E(c^{(k)})}{E'(c^{(k)})}.
\end{align}
Because $E(c) \geq 0$ then any existing root will also be a minimum implying that the multiplicity of the root is at least 2, leading to the following parameter estimation algorithm for a single parameter $c$:
\begin{align}
    c^{(k+1)} = c^{(k)} - \frac{\|I_h[v-u] \|^2_H}{\langle I_h[v-u], w \rangle_H} \bigg|_{t_k}.
    \label{eq:new1d}
\end{align}

If we consider a more specific system of the form,
\begin{align}
    \Dot{u} + \gamma Lu
    + F(u) &= 0 \notag\\
    \Dot{v} + c Lv + F(v)
     + \mu I_h (v-u) &= 0,
\end{align}
and use the asymptotic approximation \eqref{eq:asymptotic1} for the sensitivity $w$, then the  parameter update is given by
\begin{align}
    c^{(k+1)}
    = c^{(k)} + \mu
    \frac{\|I_h(v-u)\|_H^2}
    {\langle I_h (v-u) , Lv \rangle_H}
    \bigg|_{t_k},
\end{align}
which is precisely the algorithm proposed in \cite{CHL_2020} and further justified in \cite{martinez2022convergence,biswas2023determining} for the 2D Navier-Stokes equations.

When the set of unknown parameters $c$ is higher dimensional then the direct extension of \eqref{eq:new1d} is modified to become:
\begin{align}
    c^{(k+1)}_i
    = c_i^{(k)} - 
    \frac{\|I_h(v-u)\|^2_{H}}
    { \| \langle I_h(v-u), \mathbf{w} \rangle_{H} \|^2_{\mathbb{R}^n}}
    \cdot
    \langle I_h(v-u), w_i\rangle_H \bigg|_{t_k}
    .
    \label{eq:newMP}
\end{align}
We can estimate the sensitivity $w_i$ either by direct computation or by using the perturbative approximation \eqref{eq:asymptotic1}.
Note that rootfinding for functions $f:\mathbb{R}^n \to \mathbb{R}$ will not always have a unique solution.
The local linearized system will be underdetermined.
In \eqref{eq:newMP} we are using steepest descent to determine the direction in the parameter space.
More sophisticated root-finding approaches are possible however we shall not discuss them here, as we instead show that a more general way to approach the problem is as an optimization problem rather than a root-finding one. 


\subsection{Optimization methods}
We would like to make use of the Newton-Raphson method \cite{bertsekas1997nonlinear}
\begin{equation}
c^{(k+1)} = c^{(k)} - (D^2 E)^{-1}\nabla E(c),
\end{equation}
which requires the Hessian or matrix of second order derivatives, i.e.
\begin{align}
    (D^2 E)_{ij}
    = \langle w_i, I_h w_j \rangle_H \bigg|_{T}
    + \left\langle I_h (v-u), \frac{\partial^2}{\partial c_i \partial c_j}
    v \right\rangle_H
\end{align}
Implementing this as a parameter estimation algorithm requires knowledge of all of the second derivatives with respect to the approximate parameters.
However, if we make the assumption that either $v - u$  or $\frac{\partial^2}{\partial c_i \partial c_j} v$ is small, then we have the approximation of the Hessian,
\begin{align}
    (D^2 E)_{ij}
    \approx \langle w_i, I_h w_j \rangle_H \bigg|_{T}
\end{align}
This leads to Gauss Newton \cite{humpherys2020foundations} algorithm, which in our case is given by 
\begin{align}
        c_i^{(k+1)}
    = c_i^{(k)} - 
    \sum_j
    \left(
    \langle I_hw_i,I_hw_j \rangle_H
    \right)^{-1}
     \langle v-u, I_h w_j \rangle_H
     \bigg|_{t_k}.
\end{align}
The Levenberg-Marquardt algorithm is a modification of the Gauss-Newton algorithm.
It is generally more robust than the Gauss-Newton algorithm \cite{marquardt1963algorithm}; for many cases, it will converge to the correct answer even when the initial guess is far off.
It may also help in points where the matrix $\langle w_i, w_j \rangle$ is not invertible.
It is given by the following:
\begin{align}
    c_i^{(k+1)}
    = c_i^{(k)} - 
    \sum_j
    \left(
    \langle I_hw_i,I_hw_j \rangle_H
    + \lambda \delta_{ij}
    \right)^{-1}
     \langle v-u, I_h w_j \rangle_H
     \bigg|_{t_k}.
\end{align}
The inverse given above will always be defined so long as $-\lambda$ is not an eigenvalue of $\langle I_hw_i,I_hw_j \rangle_H$.
The matrix $\langle I_hw_i,I_hw_j \rangle_H$ will be positive semi-definite and so this algorithm will always be defined. 
We can estimate the sensitivity $w_i$ either by direct computation or by using the perturbative approximation \eqref{eq:asymptotic1}.


\section{Examples}
To demonstrate the utility of the OTF parameter estimation algorithm we will consider several examples of varying complexity.  Each of these examples will compare the OTF approach where the asymptotic approximation is utilized, with the DS approach where the sensitivity equations themselves are simulated.  In addition, we report on the influence and effect of several other hyper-parameters that arise in the selected algorithm.

\subsection{Lorenz '63 System}

We first consider the Lorenz '63 system originally proposed in \cite{lorenz1963deterministic}, here rewritten to match the notation established in this article as:
\begin{align}
    \Dot{u}_1 &=  -\gamma_1 (u_1-u_2) 
    \notag\\
    \Dot{u}_2 &=
    u_1(\gamma_2 - u_3) - u_2 \notag\\
    \Dot{u}_3 &=
    u_1 u_2 - \gamma_3 u_3.
\end{align}
Classical representations of this system have the following notational substitutions: $u_1 \rightarrow x,~u_2 \rightarrow y,~u_3\rightarrow z, \gamma_1 \rightarrow \sigma,~\gamma_2\rightarrow \rho,$ and $\gamma_3\rightarrow \beta$.
The corresponding nudged system is 
\begin{align}
    \Dot{v}_1 &= -c_1
    (v_1-v_2) 
    - \mu (v_1 -u_1)
    \notag\\
    \Dot{v}_2 &=
    v_1(c_2 - v_3) - v_2
    - \mu (v_2 -u_2)
    \notag\\
    \Dot{v}_3 &=
    v_1 v_2 - c_3 v_3
    - \mu (v_3 -u_3).
\end{align}

As noted in \cite{carlson2021dynamically}, the parameter update appears to converge only if the number of parameters being estimated is less than or equal to the rank of the observation operator. 

In practice, this means that if we want to estimate all three parameters in the Lorenz system, we have to nudge all of the variables. 
To focus on the parameter estimation rather than the convergence properties of the CDA algorithm for this system, we suppose that the entire state is observed.
Nudging with some subset of the state was thoroughly explored in \cite{blocher2018data,du2021analysis,carlson2021dynamically}
Indeed, \cite{carlson2021dynamically} looked at parameter estimation using a version of the algorithm proposed in \cite{CHL_2020} to do parameter estimation on various subsets of the full parameters.

The sensitivity equations for the Lorenz system are given by the system of 9 differential equations for  $w_{ij} = \frac{\partial v_i}{\partial c_j}$, i.e.
\begin{align}
    \Dot{w}_{1j} &= 
    - \delta_{1j}
    (v_1-v_2)
    - c_1 (w_{1j} - w_{2j})
    - \mu_1 w_{1j},
    \notag\\
    \Dot{w}_{2j} &=
    w_{1j}(c_2 - v_3) 
    +v_1(\delta_{2j} - w_{3j}) 
    - w_{2j}
    - \mu_2 w_{2j},
    \notag\\
    \Dot{w}_{3j} &=
    w_{1j}v_2
    + v_1 w_{2j}
    - \delta_{3j} v_3
    - c_3 w_{3j}
    - \mu_3 w_{3j}.
\end{align}
Note that because we have the same number of variables as parameters we can represent $w_{ij}$ as a square matrix, $W = [w_{ij}]$.
In this case we get that the gradient of the error $ E(c) = \frac{1}{2} \|\mathbf{v}-\mathbf{u}\|^2_{\mathbb{R}^3} \big|_{t=T}$ is
\begin{align}
    \nabla E
    = W^T (\mathbf{v} - \mathbf{u}) \big|_{t=T}.
\end{align}

In each of the optimization routines discussed below, we take the relaxation parameter to be $\mu = 100$ 
unless otherwise specified, and the time between parameter updates as $\Delta t = 0.5$.
We take the true parameter values as the canonical $\gamma = (10,28,8/3)$, and start with initial conditions $\mathbf{u} = (0,1,-1)$ and $\mathbf{v} = \mathbf{0}$.
We start with initial parameter values $c = \frac{1}{2} \gamma$.

\begin{itemize}
\item The gradient descent algorithm in this setting is given by
\begin{align}\label{eq:L63_GD}
    c^{k+1} = c^{k} -  r W^T (\mathbf{v} - \mathbf{u}) \big|_{t =t_k},
\end{align}
where the learning rate $r$ is selected to provide the optimal rate of convergence.  Adjusting the learning rate $r$ leads to different rates of convergence up to $r\approx 40$ whereupon$\mathbf{v}$ no longer converges to $\mathbf{u}$. 
For the Lorenz '63 model considered here, we set $r=30$.

\item  Newton's method for the Lorenz '63 system gives the following update: 
\begin{align}\label{eq:L63_Newton}
    c^{(k+1)}
    = 
    c^{(k)}
    -
    \frac{\|\mathbf{v}- \mathbf{u}\|^2}
    {\| W^T(\mathbf{v}-\mathbf{u})\|^2 }
    W^T (\mathbf{v}-\mathbf{u}).
\end{align}

\item The Levenberg-Marquardt algorithm is given by
\begin{align}\label{eq:L63_LM}
    c^{(k+1)}
    &= 
    c^{(k)}
    -
    (W^TW + \lambda I)^{-1}
    W^T (\mathbf{v}-\mathbf{u}).
\end{align}
We select $\lambda = 10^{-6}$ for this system. 
\end{itemize}

\begin{figure}
    \centering
    \includegraphics[width=0.9\textwidth]
    {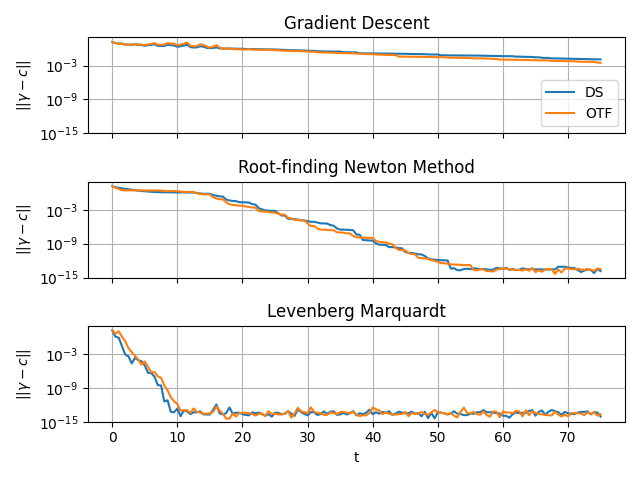}
    \caption{Comparison of DS and OTF methods for the Lorenz '63 model.  Gradient descent is implemented with a learning rate of $r=30$, and the smoothing parameter in Levenberg-Marquardt is $\lambda = 10^{-6}$ for both the OTF and DS methods.  The default value of $\mu = 100$ is also used for both, with parameter updates occurring at every $\Delta t= 0.5$ time units.
    }
    \label{fig:OTFvsDS_lorenz}
\end{figure}

For the Lorenz '63 system, the asymptotics in Section \ref{sec:asym} lead to the following approximation to the parameter sensitivities:
\begin{align}
    W^{OTF} = 
    \begin{bmatrix}
        -\frac{1}{\mu_1}(v_1-v_2) & 0 & 0 \\
        0 & \frac{1}{\mu_2}v_1 & 0 \\
        0 & 0 & - \frac{1}{\mu_3} v_3
    \end{bmatrix}.
\end{align}
Replacing $W\rightarrow \tilde{W}$ in \eqref{eq:L63_GD}, \eqref{eq:L63_Newton}, or \eqref{eq:L63_LM} produces the OTF version of each of these algorithms.  In Figure \ref{fig:OTFvsDS_lorenz}, we can see how each of these methods compares (with the same default learning rate $r$ and smoothing parameter $\lambda$ as given above).  Note that there is little noticeable difference in the convergence rate or final achieved error for OTF versus DS, despite the significant computational savings in the OTF method.

\begin{figure}
    \centering
    \includegraphics[width=0.9\textwidth]
    {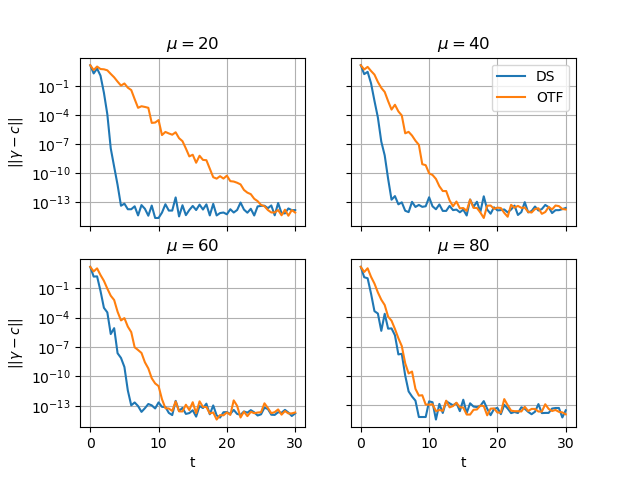}
    \caption{DS vs. OTF Levenberg-Marquardt method for various values of $\mu$. All other parameters are the same as those in Figure \ref{fig:OTFvsDS_lorenz}.  Note the similar behavior for larger values of $\mu$ in agreement with perturbative results. }
    \label{fig:LM_mu}
\end{figure}
As the OTF method is an asymptotic approximation as $\mu\rightarrow \infty$, we anticipate that the two methods will agree for larger values of $\mu$. 
In Figure \ref{fig:LM_mu} we compare the convergence of the parameter error for Levenberg-Marquadt DS versus OTF at different values of the nudging parameter $\mu$ (all other parameters are the same as in Figure \ref{fig:OTFvsDS_lorenz}).  Note that even for moderate values of $\mu$, the OTF method achieves the same level of accuracy nearly as quickly as DS.  In fact, even for $\mu=20$, OTF still achieves the same level of accuracy, the convergence rate is just significantly slower.

\subsection{Two Layer Lorenz '96 Model}
We next consider the two layer Lorenz '96 model originally derived to model atmospheric quantities along a single circle of latitude \cite{lorenz1996predictability}, and modified to include multiple scales in \cite{chen2017beating,chen2018conditional}. The original Lorenz '96 model has been used as a testbed for various data assimilation techniques since its inception (see \cite{anderson2001ensemble,ott2004local,wilks2005effects,arnold2013stochastic} for example) including modern machine learning techniques (see \cite{chattopadhyay2020data} for one example).  More recently, a modification of the parameter estimation algorithm originally proposed in \cite{pachev2022concurrent} was applied to the two-layer system in \cite{martinez2024relaxation} to indicate the utility of parameter recovery/estimation when a large number of parameters were to be estimated.  

The Lorenz '96 system provides a high, yet finite dimensional dynamical system with two different scales that are coupled together.
$u^l$ captures the large scale behavior of the system, with $k \in \{1,..,I\}$ indexing the position along the latitude circle.
Each $u^l_k$ is coupled to $J$ small scale variables $u^s$, which represent the small scale behavior not accounted for by the coarse grained ($u^l$ only) system.
\begin{align}
    \Dot{u}_k^l &= u_{k+1}^l
    (u_{k-1}^l - u_{k+2}^l)
    + \gamma_1\sum_{j=1}^J u_{kj}^s u_{k}^l 
    - \gamma_2 u_k^l + F \notag\\
    \Dot{u}_{kj}^s,
    &= - d_j u_{kj}^s - \gamma_1 (u_k^l)^2.
\end{align}
This is coupled to the data assimilated system:
\begin{align}
    \Dot{v}_k^l &= v_{k+1}^l
    (v_{k-1}^l - v_{k+2}^l)
    + c_1\sum_{j=1}^J v_{kj}^s v_{k}^l 
    - c_2 v_k^l + F 
    - \mu ( v^l_k -u^l_k),
    \notag\\
    \Dot{v}_{kj}^s
    &= - d_j v_{kj}^s - c_1 (v_k^l)^2.
\end{align}
We will take the large scale variables, i.e. the $u^l$ as the observable states of the system so that the observation operator is a projection from the entire state $(u^l,u^s)$ onto $u^l$ only.
Since the parameters $d_j$ only appear in the small scale equations, they only directly interact with the unobservable variables, and hence we can not anticipate recovering their true value with these relevant observations.
Following \cite{chen2018conditional}, we take these parameters fixed as $(d_j) = (0.2,0.5,1,2,5)$.
We take the other true parameter values as $\gamma = (0.01,0.5)$, with starting approximate parameters as $c = \frac{1}{2}\gamma$.
We also select $J=5$ small scale variables for each large scale variable, and and $I = 40$ large scale variables in total.
For this problem we will use a nudging parameter of $\mu = 50$ which is sufficient for the system with correct parameters to converge within about 175 time units.
Our initial conditions for $u$ are randomly drawn from a normal distribution, and we take 0 initial conditions for $v$.

\begin{figure}
    \centering
    \includegraphics[width=0.9\textwidth]{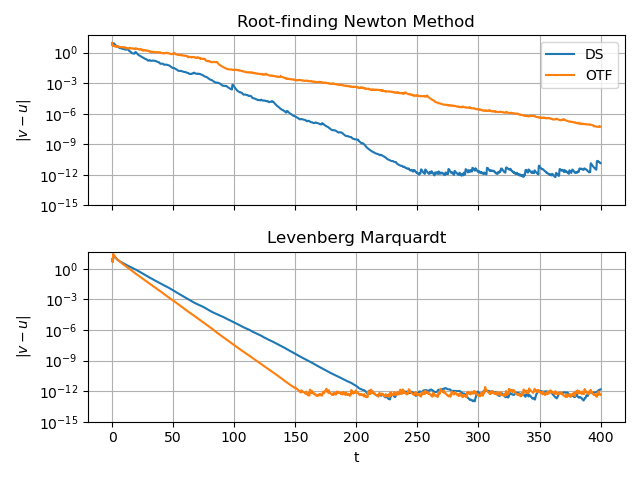}
    \caption{Comparison between the DS and OTF methods for the Lorenz '96 model.  Note that in this case the asymptotic approximation utilized in the OTF model appears to outperforming the DS approach.}
    \label{fig:L96_compare}
\end{figure}

The sensitivity equations for Lorenz '96 are given by:
\begin{align}
    \Dot{w}_{k,i}^l &= w_{k+1,i}^l
    (v_{k-1}^l - v_{k+2}^l)
    + v_{k+1}^l
    (w_{k-1,i}^l - w_{k+2,i}^l) 
    + \delta_{1i}\sum_{j=1}^J v_{kj}^s v_{k}^l,
    \notag\\
    &+ c_1\sum_{j=1}^J w_{kj,i}^s v_{k}^l
    + c_1\sum_{j=1}^J v_{kj}^s w_{k,i}^l
    - \delta_{2i} v_k^l  
    - c_2 w_{k,i}^l  
    - \mu w^l_{k,i}
    \notag\\
    \Dot{w}_{kj,i}^s
    &= - d_j w_{kj,i}^s 
    - \delta_{1,i} (v_k^l)^2
    - 2 c_1 v_k^l w_{k,i}^l.
\end{align}
For this model, we consider only the implementation of Newton's method and the Levenberg-Marquardt algorithm.
Following the asymptotics developed above, we identify the following approximations for the sensitivities as $\mu\to \infty$:
\begin{align}
    w_{k,1}^{l,OTF} = \frac{1}{\mu}\sum_{j=1}^J v_{kj}^s v_{k}^l,
    \notag\\
    w_{k,2}^{l,OTF} = -\frac{1}{\mu} v^l_k.
\end{align}
This is compared to the DS method in Figure \ref{fig:L96_compare} where the default values of $\mu = 50$, $\delta t = 0.5$ and $\lambda = 10^{-6}$ in Levenberg-Marquadt are used.
Interestingly for the Lorenz '96 model, these parameter estimation algorithms occasionally will fail to converge to the true value, and instead converge to $(c_1,c_2) = (-\gamma_1,\gamma_2)$ which appears to be an underlying feature of the symmetries inherent to the system itself.  This can of course be avoided with careful selection of the hyper-parameters ($\mu$ and/or $\Delta t$) or the initial guess for $c_1$ and $c_2$, but also illustrates the ill-posedness of the underlying problem.


\subsection{Kuramoto-Sivashinsky Equation}
Having shown the applicability of the OTF approach to two finite dimensional systems, we now turn to the infinite dimensional case to study how these methods can be applied to partial differential equations as well.
We consider the one (spatial) dimensional Kuramoto-Sivashinsky Equation (KSE)
\begin{align}\label{eq:KSE_truth}
    \Dot{u}
    + \gamma_1 u''
e    + \gamma_2 u u'
    +\gamma_3 u^{(4)}
     = 0
\end{align}
here we have $u' = \frac{\partial u}{\partial x}$, and the higher order derivatives are defined similarly.
The KSE system is used to model a variety of physical systems. Generally, it models systems  far from equillibrium such as turbulent behavior of a single flame, instabilities in reaction diffusion systems, and the flow of plasmas \cite{kuramoto1978diffusion,ashinsky1988nonlinear}.
It is mathematically interesting as a model because it is an example of a chaotic PDE (see \cite{nicolaenko1985some,hyman1986kuramoto,collet1993global,misbah1994secondary,cvitanovic2010state} for example) with only one spatial dimension, and as such, is relatively computationally inexpensive to simulate. Parameter estimation in the KSE system was studied under a similar framework in \cite{pachev2022concurrent}, and a similar question was posed for KSE in \cite{mojgani2022discovery} although a very different approach to the parameter estimation problem was taken there.

The data assimilated/nudged system for KSE is given by
\begin{align}\label{eq:KSE_nudge}
    \Dot{v}
    + c_1 v''
    + c_2 v v'
    + c_3 v^{(4)}
    + \mu I_h(v-u)
    = 0.
\end{align}
In this example we will use the observation operator $I_h$ as a projection onto the lowest $N = 1/h$ Fourier modes of the solution.
In all that follows, we will use inner products defined via the $L^2$ norm.
In everything below we will let $N = 32$ (the full simulations run with a resolution of 1024 grid points).
Because we have three parameters to fit, DS methods would greatly increase the computational complexity of the algorithm, resulting in 4 PDE solves per time step instead of a single solve for a first order DS method.
We first fit only $c_1$ where DS methods do not greatly increase the computational cost.
This allows us to compare the DS methods with the OTF methods for a single parameter recovery.  Following this investigation into single parameter estimation, we simulate multi-parameter estimation using only OTF methods which greatly reduces the computational cost.

\subsubsection{Single Parameter Estimation}
For estimation of the single parameter $c_1$ we need the evolution equation for the first order sensitivity equation: 
\begin{align}
    &\Dot{w}
    + v''
    + c_1 w''
    + c_2 (w v)' 
    + c_3 w^{(4)}
    + \mu P_N w = 0,
\end{align}
We take the true parameter values $\gamma = (1,1,1)$, and the hyper parameters $\Delta t = 0.5$ and $\mu = 25$.
We take initial conditions
\begin{multline}
    u|_{t=0}
   = \sin(6 \pi x/L ) + 0.1 \cos(\pi x/L )
          - 0.2 \sin(3\pi x/L ) \\
          + 0.05\cos(15 \pi  x/L )
          + 0.7 \sin(18 \pi x/L )
          - \cos(13 \pi x/L ).
\end{multline}
where our domain is $[0,L] = [0,100]$, and we take the initial conditions of $v$ to be identically zero.  Periodic spatial boundaries are implemented throughout this comparison and all that follows.

\begin{figure}
    \centering
    \includegraphics[width=0.85\textwidth]
    {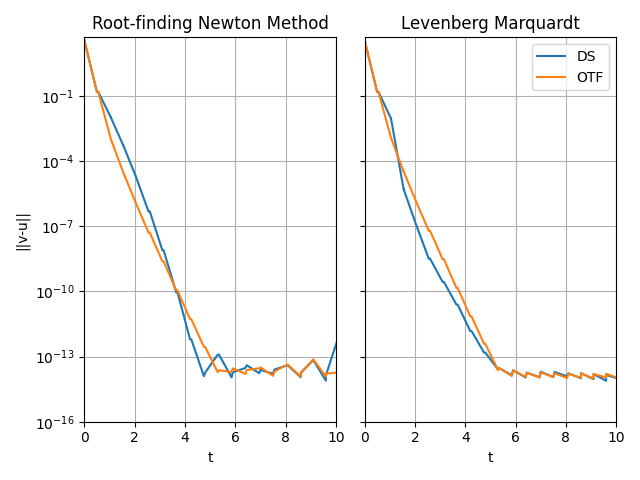}
    \caption{Comparison of OTF and DS methods for estimating a single parameter for KSE with a nudging parameter $\mu = 25$ where $N=32$ modes are observed, and the parameter is updated every $\Delta t = 0.5$ time units. 
    }
    \label{fig:KSE_single_DSvsOTF}
\end{figure}

Directly simulating the sensitivity equations and applying Newton's root-finding method and Levenberg-Marquardt in this setting works quite well as indicated in Figure \ref{fig:KSE_single_DSvsOTF} (with initial guess $c_1 = \frac{1}{2} \gamma_1$).  
Note that Newton's root solving method marginally converges faster in this setting, although the differences are quite slight.  
We also find that the optimization-based routines such as Levenberg-Marquardt and gradient descent are more sensitive to the initial choice of $c_1$ then Newton's root-finding method.  
In fact, with an initial value of $c_1 = 4.5$ the optimization based methods converge to $c_1 \approx 4.25$ whereas Newton's root finding method still manages to recover the true value.  


We compare the DS and OTF methods for single parameter estimation in KSE in Figure \ref{fig:KSE_single_DSvsOTF}.  Note that even for $\mu=25$ as used here, OTF follows the DS estimate at a fraction of the cost.

\subsubsection{OTF Methods For Multiparameter Estimation}
To perform multiparameter estimation we need the full sensitivity equations for this system which are given by:
\begin{align}
    \Dot{w}_i 
    +
    \delta_{1i} v''
    + \delta_{2i} v v'
    + \delta_{3i} v^{(4)}
    + c_1 w_i''
    + c_2 w_i v'
    +c_2 v w_i'
    +c_3 w_i^{(4)}  + \mu P_N w_i
    =0.
\end{align}
The asymptotic approximation for each of these terms is given by:
\begin{align}
    \begin{bmatrix}
        w_1^{OTF} \\
        w_2^{OTF} \\
        w_3^{OTF}
    \end{bmatrix}
    \sim 
    -\frac{1}{\mu}
    \begin{bmatrix}
    v'' \\
    v v' \\
    v^{(4)}
    \end{bmatrix}.
\end{align}
This gives the Levenberg-Marquardt approximation to the Hessian:
\begin{align}
    D^2 E_{LM}^{OTF}
    &
    =
    \frac{1}{\mu}
    \begin{bmatrix}
    \|I_hv''\|^2
    + \lambda
    &
    \langle v'', I_h[vv'] \rangle
    & 
    \langle v'',I_h [v^{(4)}]\rangle
    \\
    .
    & 
    \|I_h[vv']\|^2
    + \lambda
    & 
    \langle v^{(4)}, I_h[vv'] \rangle
    \\
    .
    &.
    &
    \|I_h[v^{(4)}]\|^2
    +\lambda
    \end{bmatrix}.
\end{align}

\begin{figure}
    \centering
    \includegraphics[width=0.95\textwidth]
    {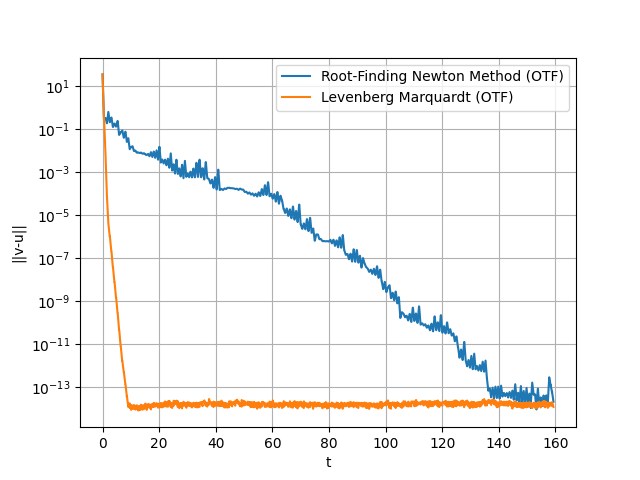}
    \caption{Comparison of algorithms for estimating  all three parameters $(c_1,c_2,c_3)$ simultaneously for KSE.  In this setting the initial parameter values are double the true values, the time between updates is $\Delta t = 0.5$, and the nudging coefficient is $\mu = 25$ with the observation operator $I_h$ defined as a projection onto the lowest 32 Fourier modes of the solution.}
    \label{fig:NvsLM_KSE}
\end{figure}

We take $c^{(0)} = (2,2,2)$ as the initial guess where the true parameter values are $\gamma = (1,1,1)$, and the hyper parameters $\Delta t = 0.5$ and $\mu = 25$ with $I_h$ being defined as the projection onto the lowest 32 Fourier modes.  The results are summarized in Figure \ref{fig:NvsLM_KSE} for the Newton root-finding method and Levenberg Marquardt  (with $\lambda = 10^{-6}$) both via the asymptotic OTF approximation.  Note that although both methods converge, Levenberg Marquradt does so much faster, likely because root finding in high dimensions is a notoriously difficult problem.


\subsubsection{Parameter Identification in the absence of a `true' value} \label{sec:optvsroot}
We now investigate the situation where the underlying model assumption is incorrect, that is if the ansatz we propose for the underlying observed data is incorrect, we want to consider how parameter estimation will perform.  Specifically, we will replace the `true' system in \eqref{eq:KSE_truth} with the following
\begin{align}\label{eq:KSE_LES}
    \Dot{u}
    + \gamma_1 u''
    + \gamma_2 u u'
    +\gamma_3 u^{(4)}
    - \epsilon u^{(6)}
     = 0.
\end{align}
We will let $\epsilon = 10^{-3}$ be fixed, which introduces a perturbation of the original system \eqref{eq:KSE_truth}.  Here we will take \eqref{eq:KSE_nudge} as the data assimilated version of this system and identify the corresponding parameters $(c_1,c_2,c_3)$, i.e. we are intentionally trying to assimilate data into a model that is inherently incorrect.  In this setting, the additional sixth order term can be thought of as a representation of some process for which the model doesn't account, i.e. a significant process for which the original model developer was unaware.  This term may also be a stand-in for issues related to finite numerical resolution of the solution, i.e. if \eqref{eq:KSE_truth} is not adequately resolved, one way of representing this under-resolved solution is to consider additional dissipation such as the higher order term introduced here.  This is in line with the large eddy simulation (or subgrid scale) models that are in use for fluid dynamics (see \cite{jakhar2024learning,guan2024online} for example)

\begin{figure}
    \centering
    \includegraphics[width=0.95\textwidth]
    {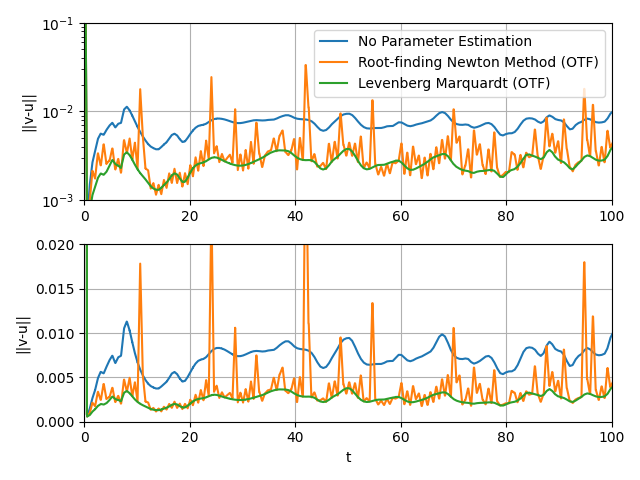}
    \caption{
    Parameter estiimation for the modified KSE (see \eqref{eq:KSE_LES}) with $\epsilon = 10^{-3}$ .  The blue curve represents the effects of performing nudging on this system without parameter updates applied.  Note that although the Newton root finding method roughly agrees with Levenberg-Marquadt, the root finding technique employed here is far noisier both for the parameter updates, and the state estimate.}
    \label{fig:KSE_optim}
\end{figure}

For this setting, we will initialize the parameter $c_i$ at the `true' values, i.e. $c_i = \gamma_i = 1$ and see how these parameters are adjusted to account for the nonzero presence of $\epsilon$ in the generation of the observable data.
We see in Figure \ref{fig:KSE_optim} that the Levenberg-Marquardt algorithm effectively reduces the error by a factor of $3$.  On the other hand, Newton's root finding method produces oscillating estimates of the parameters, yielding an inconsistent estimate of both the underlying parameter and hence the state itself.  This is not unexpected, as $E(c_1,c_2,c_3)$ will not have a root in this case.

\section{Conclusion}


We have outlined an original framework for parameter estimation within the continuous time data assimilation approach.
Starting from some basic assumptions, motivated by numerical experiments, we have reduced the problem to a finite dimensional optimization problem.
Gradient based optimization methods applied to this setting lead to the sensitivity equations for the desired parameter, set in the data assimilation framework.
Asymptotic expansion with respect to the nudging coefficient on these sensitivity equations yields computationally efficient approximations that imitate the full sensitivity equations which yields an on-the-fly (OTF) parameter estimation algorithm.
In particular, the OTF version of Newton's root-finding method for one parameter reproduces the algorithm originally proposed in \cite{carlson2021dynamically}.  Using various optimization routines (rather than Newton's root-finding technique), we have developed a suite of parameter estimation algorithms that accurately identify the unknown parameters both for multiple unknown parameters and for parameters that appear in nonlinear operators in the underlying dynamical system, thus avoiding pitfalls that were observed in similar algorithms proposed in \cite{pachev2022concurrent} and elsewhere. This parameter estimation framework has been demonstrated to work well on 3 specific models of increasing complexity, and has demonstrated a certain level of utility even in the case where the full model is incorrectly specified.  

Future work will establish the rigorous justification of this approach, providing rigorous estimates on the asymptotic approximation which lies at the core of this article, justifying the discrete time updates of the parameter, and clarifying the role that the feedback control data assimilation method plays in this approach.  Further application of this method will be to adapt it for estimation of parameters when the observables are stochastically perturbed, extension to more complicated settings where the parameters are unknown, modifications for different norms when defining the optimization objective $E(c)$, and applications to equation discovery.

\section*{Acknowledgements}
JPW was partially supported by NSF grant DMS-2206762 and CCF-343286.  JN was partially supported by NSF grant DMS-2206762.  E.C. was partially supported  by the Department of Defense Vannevar Bush Faculty Fellowship, under ONR award N00014-22-1-2790.


\bibliographystyle{plain}


\end{document}